\documentstyle[11pt]{article}
\def \reel{ {\rm I}\!{\rm R} }
 \newcommand{\too}{\longrightarrow}
\newcommand{\om}{\omega}
\newcommand{\Om}{\Omega}
\newcommand{\na}{\nabla}

\newcommand{\al}{\alpha}
\newcommand{\la}{\lambda}
\newcommand{\be}{\beta}
\newcommand{\ga}{\gamma}

\newcommand{\G}{{\cal G}}

\newcommand{\U}{{\cal U}}

\newcommand{\de}{\delta}

\newcommand{\inj}{\hookrightarrow}

 \def \rat{ {\rm Q}\kern-.65em {}^{{}_/ }}

\newtheorem{th}{Theorem}[section]
\newtheorem{pr}{Proposition}[section]
\newtheorem{Le}{Lemma}[section]
\newtheorem{Rm}{Remark}[section]

\title{On the Riemann-Lie algebras and Riemann-Poisson Lie groups} \author{M. Boucetta}
\date{} \parindent=0cm
\begin{document}
\maketitle

{\bf Abstract.} A Riemann-Lie algebra is a Lie algebra $\G$ such
that its dual $\G^*$ carries a Riemannian metric compatible (in
the sense introduced by the author in [1]) with the canonical
linear Poisson structure  of $\G^*$.  The notion of Riemann-Lie
algebra has its origins in the study, by the author, of
Riemann-Poisson manifolds (see [2]).

In this paper, we  show that, for a Lie group $G$, its  Lie
algebra $\G$ carries a structure of  Riemann-Lie algebra iff
 $G$ carries a flat left-invariant Riemannian metric.
 We
use this characterization to construct  a huge number of
Riemann-Poisson Lie groups (a Riemann-Poisson Lie group is a
Poisson Lie group endowed with a left-invariant Riemannian metric
compatible with the Poisson structure).
\section{Introduction}
 Riemann-Lie algebras first arose in the study by
the author of  Riemann-Poisson manifolds (see [2]). A
Riemann-Poisson manifold is a Poisson manifold $(P,\pi)$ endowed
with a Riemannian metric $ <,>$  such that the couple $(\pi,<,>)$
is compatible in the sense introduced by the author in [1]. The
notion of Riemann-Lie algebra appeared when we looked for a
Riemannian metric compatible with the canonical Poisson structure
on the dual of a Lie algebra. We pointed out (see [2]) that the
dual space $\G^*$ of a Lie algebra $\G$ carries a Riemannian
metric compatible with the linear Poisson structure iff the Lie
algebra $\G$ carries a structure which we called  Riemann-Lie
algebra.
 Moreover, the isotropy Lie algebra at a point on a
Riemann-Poisson manifold is a Riemann-Lie algebra. In particular,
the dual Lie algebra of a Riemann-Poisson Lie group is a
Riemann-Lie algebra (a Riemann-Poisson Lie group is a  Poisson Lie
group endowed with a left-invariant Riemannian metric compatible
with the Poisson structure). In this paper, we will show that a
Lie algebra $\G$ carries a structure of Riemann-Lie algebra iff
$\G$ is a semi-direct product of two abelian Lie algebras. Hence,
according to a well-known result of Milnor [5], we deduce that,
for a Lie group $G$, its Lie algebra carries a structure of
Riemann-Lie algebra iff  $G$ carries a flat left-invariant
Riemannian metric. As application, we give a method to construct a
huge number of Riemann-Poisson Lie groups.  In particular, we give
many examples of bialgebras $(\G,[\;,\;],\G^*,[\;,\;]^*)$ such
that both $(\G,[\;,\;])$ and $(\G^*,[\;,\;]^*)$ are Riemann-Lie
algebras.

\section{Definitions and main results}
\subsection{Notations }
Let $G$ be a connected Lie group and  $(\G,[\;,\;])$ its Lie
algebra. For any $u\in\G$, we denote by $u^l$ (resp. $u^r$) the
left-invariant (resp. right-invariant) vector field of $G$
corresponding to $u$. We denote by $\theta$  the right-invariant
Maurer-Cartan form on $G$ given by
$$\theta(u^r)=-u,\quad u\in\G.\eqno(1)$$

Let $<;>$ be a scalar product on $\G$ i.e. a bilinear, symmetric,
non-degenerate and  positive definite form on $\G$.

Let us enumerate some mathematical objets which are naturally
associated with $<,>$:

1. an isomorphism $\#:\G^*\too\G$;

2. a scalar product $<,>^*$ on the dual space $\G^*$ by
$$<\al,\be>^*=<\#(\al),\#(\be)>\qquad\al,\be\in\G^*;\eqno(2)$$

3. a left-invariant Riemannian metric $<,>^l$ on $G$ by
$$<u^l,v^l>^l=<u,v>\qquad u,v\in\G;\eqno(3)$$

4. a left-invariant contravariant Riemannian metric $<,>^{*l}$ on
$G$ by
$$<\al,\be>^{*l}_g=<T_e^*L_g(\al),T_e^*L_g(\be)>^*\eqno(4)$$where
$\al,\be\in\Om^1(G)$ and $L_g$ is the left translation of $G$ by
$g$.

The infinitesimal  Levi-Civita connection associated with
$(\G,[\;,\;],<,>)$ is the  bilinear map $A:{\cal G}\times{\cal
G}\too{\cal G}$ given by
$$2<A_uv,w>=<[u,v],w>+<[w,u],v>+<[w,v],u>,
\quad u,v,w\in\G.\eqno(5)$$ Note that $A$ is the unique bilinear
map from $\G\times\G$ to $\G$ which verifies:

1. $A_uv-A_vu=[u,v]$;

2. for any $u\in\G$, $A_u:\G\too\G$ is skew-adjoint i.e.
$$<A_uv,w>+<v,A_uw>=0,\quad v,w\in\G.$$
 The terminology used here can be justified by the fact that the
 Levi-Civita connection $\na$ associated
with $(G,<,>^l)$ is given by
$$\na_{u^l}v^l=(A_uv)^l, u,v\in\G.\eqno(6)$$
We will introduce now a Lie subalgebra of $\G$ which will play a
crucial role in this paper.

For any $u\in\G$, we denote by $ad_u:\G\too\G$ the endomorphism
given by $ad_u(v)=[u,v]$, and  by $ad^t_u:\G\too\G$ its adjoint
given by $$<ad^t_u(v),w>=<v,ad_u(w)>\qquad v,w\in\G.$$ The space
$$S_{<,>}=\{u\in\G;ad_u+ad^t_u=0\}\eqno(7)$$ is obviously a
subalgebra of $\G$. We call $S_{<,>}$ {\it the orthogonal
subalgebra associated with $(\G,[\;,\;],<,>)$}.

\begin{Rm}
 The scalar product $<,>$ is bi-invariant if and only if $\G=S_{<,>}$.
In this case $\G$ is the product of an abelian Lie algebra and a
semi-simple and compact Lie algebra (see [5]). The general case
where $<,>$ is not positive definite has been studied by A. Medina
and P. Revoy  in [4] and they called a Lie algebra $\G$ with an
inner product $<,>$ such that $\G=S_{<,>}$ an orthogonal Lie
algebra which justifies the terminology used here.$\Box$\end{Rm}

Let $({\cal G},[\;,\;],<,>)$ be a Lie algebra endowed with a
scalar product.

The triple $({\cal G},[\;,\;],<,>)$ is called {\it a Riemann-Lie
algebra} if
$$[A_uv,w]+[u,A_wv]=0\eqno(8)$$for all $u,v,w\in\cal G$ and where
 $A$ is
the infinitesimal Levi-Civita connection associated to $({\cal
G},[\;,\;],<,>)$.

From the relation $A_uv-A_vu=[u,v]$ and the Jacobi identity, we
 deduce that (8) is equivalent to
$$[u,[v,w]]=[A_uv,w]+[v,A_uw]\eqno(9)$$for any $u,v,w\in\cal G$.
We refer the reader to [2] for the origins of this definition.

Briefly, if $({\cal G},[\;,\;],<,>)$ is a Lie algebra endowed with
a scalar product. The scalar product $<,>$ defines naturally a
contravariant Riemannian metric on the Poisson manifold $\G^*$
which we denote also by $<,>$.  If we denote by $\pi^l$ the
canonical Poisson tensor on $\G^*$, $(\G^*,\pi^l,<,>)$ is a
Riemann-Poisson manifold iff the triple $({\cal G},[\;,\;],<,>)$
is  a Riemann-Lie algebra.

\subsection{Characterization of Riemann-Lie algebras}

With the notations and the definitions above, we can state now the
main result of this paper.
\begin{th} Let $G$ be a Lie group, $(\G,[\;,\;])$ its Lie algebra
and $<,>$ a scalar product on $\G$. Then, the following assertions
are equivalent:

1) $(\G,[\;,\;],<,>)$ is a Riemann-Lie algebra.

2) $(\G^*,\pi^l,<,>)$ is a Riemann-Poisson manifold ($\pi^l$ is
the canonical Poisson tensor on $\G^*$ and $<,>$ is considered as
a contravariant metric on $\G^*$).

3) The 2-form $d\theta\in\Om^2(G,\G)$ is parallel with respect the
Levi-Civita connection $\na$ i.e. $\nabla d\theta=0.$

4) $(G,<,>^l)$ is a flat Riemannian manifold.

5) The orthogonal subalgebra $S_{<,>}$ of $(\G,[\;,\;],<,>)$ is
abelian and $\G$ split as an orthogonal direct sum
$S_{<,>}\oplus\U$ where  $\U$ is a commutative ideal.
\end{th}

The equivalence ``$1)\Leftrightarrow2)$'' of this theorem was
proven in [2] and the equivalence  ``$4)\Leftrightarrow5)$'' was
proven by Milnor in [5]. We will prove the equivalence
``$1)\Leftrightarrow3)$'' and the equivalence
``$1)\Leftrightarrow5)$''.

The equivalence ``$1)\Leftrightarrow3)$'' is a direct consequence
of the following formula which it is easy to verify:
$$\na
d\theta(u^l,v^l,w^l)_g=Ad_g\left([u,[v,w]]-[A_uv,w]-
[v,A_uw]\right),\quad u,v,w\in\G,g\in G.\eqno(10)$$

If $G$ is compact, connected and non-abelian, the condition $\na
d\theta=0$ implies that $d\theta$ is harmonic and defines,
according to Hodge Theorem, a non-vanishing class in
$H^2(G,\reel)$.   If $G$ is compact, connected and semi-simple Lie
group, then $H^2(G,\reel)=0$ and hence we get the following lemma
which will be used in the proof of the equivalence
``$1)\Leftrightarrow5)$'' in Section 3.
\begin{Le}Let $G$ be a compact, connected and semi-simple Lie group.
Then the Lie algebra of $G$ does not admit any structure of
Riemann-Lie algebra.\end{Le}

A proof of the equivalence ``$1)\Leftrightarrow5)$'' will be given
in Section 3.
\subsection{Examples of Riemann-Poisson Lie groups}
This subsection is devoted to the construction, using  Theorem
2.1, of many structures of Riemann-Poisson Lie groups. A
Riemann-Poisson Lie group is a Poisson Lie group with a
left-invariant Riemannian metric compatible with the Poisson
structure (see [2]).

 We refer the
reader to [6] for background material on Poisson Lie groups.

Let $G$ be a Poisson Lie group with Lie algebra $\cal G$ and $\pi$
the Poisson tensor on $G$. Pulling $\pi$ back to the identity
element $e$ of $G$ by left translations, we get a map
$\pi_l:G\too{\cal G}\wedge{\cal G}$
 defined by $\pi_l(g)=(L_{g^{-1}})_*\pi(g)$
 where $(L_g)_*$
denotes the tangent map of the left  translation of $G$ by $g$.
Let
$$d_e\pi:{\cal G}\too{\cal G}\times{\cal G}$$
be the intrinsic derivative of
$\pi$ at $e$ given by
$$v\mapsto L_X\pi(e),$$ where $X$ can be any vector field on
$G$ with $X(e)=v.$

The dual map of $d_e\pi$ $$[\;,\;]_e:{\cal G}^*\times{\cal
G}^*\too{\cal G}^*$$ is exactly the Lie bracket on ${\cal G}^*$
obtained by linearizing the Poisson structure at $e$. The Lie
algebra $(\G^*,[\;,\;]_e)$ is called the dual Lie algebra
associated with the Poisson Lie group $(G,\pi)$.

 We consider now  a scalar product
$<,>^*$  on $\G^*$. We denote by $<,>^{*l}$ the left-invariant
contravariant Riemannian metric on $G$ given by $(4)$.

We have shown (cf. [2] Lemma 5.1) that $(G,\pi,<,>^{*l})$ is a
Riemann-Poisson Lie group iff, for any $\al,\be,\ga\in\G^*$ and
for any $g\in G$,
$$[Ad^*_g\left(A^*_\al\ga+ad^*_{\pi_l(g)(\al)}\ga\right),Ad^*_g(\be)]_e
+[Ad^*_g(\al),
Ad^*_g\left(A^*_\be\ga+ad^*_{\pi_l(g)(\be)}\ga\right)]_e=0,
\eqno(11)$$where $A^*:\G^*\times\G^*\too\G^*$ is the infinitesimal
Levi-Civita connection associated to $(\G^*,[\;,\;]_e,<,>^*)$ and
where $\pi_l(g)$ also denotes the linear map from $\G^*$ to $\G$
induced by $\pi_l(g)\in\G\wedge\G$.

This complicated equation can be simplified enormously in the case
where the Poisson tensor arises from a solution of the classical
Yang-Baxter equation. However, we need to work more in order to
give this simplification.

Let $G$ be a Lie group and let $r\in\G\wedge\G$. We will also
denote by $r:\G^*\too\G$ the linear map induced by $r$. Define a
bivector $\pi$ on $G$ by
$$\pi(g)=(L_g)_*r-(R_g)_*r\qquad g\in G.$$
$(G,\pi)$ is a Poisson Lie group if and only if the element
$[r,r]\in\G\wedge\G\wedge\G$ defined by
$$[r,r](\al,\be,\ga)=\al([r(\be),r(\ga)])+\be([r(\ga),r(\al)])+
\ga([r(\al),r(\be)])\eqno(12)$$ is   $ad$-invariant. In
particular, when $r$ satisfies the Yang-Baxter equation
$$[r,r]=0,\eqno(Y-B)$$it defines a Poisson Lie group structure on
$\G$ and, in this case, the bracket of the dual Lie algebra $\G^*$
is given by
$$[\al,\be]_e=ad^*_{r(\be)}\al-ad^*_{r(\al)}\be,\quad\al,\be\in\G^*.
\eqno(13)$$ We will denote by $[\;\;]_r$ this bracket.

We will give now another description of the solutions of the
Yang-Baxter equation which will be useful latter.

We observe that to give $r\in {\cal G}\wedge{\cal G}$ is
equivalent to give a vectorial subspace $S_r\subset\cal G$ and a
non-degenerate 2-form $\om_r\in\wedge^2S^*_r$.

Indeed, for $r\in {\cal G}\wedge{\cal G}$, we put $S_r=Imr$ and
$\om_r(u,v)=r(r^{-1}(u),r^{-1}(v))$ where $u,v\in S_r$ and
$r^{-1}(u)$ is any antecedent of $u$ by $r$.

Conversely, let $(S,\om)$ be  a vectorial subspace of $\cal G$
with a non-degenerate 2-form. The 2-form $\om$ defines an
isomorphism $\om^b:S\too S^*$ by $\om^b(u)=\om(u,.)$, we denote by
$\#:S^*\too S$ its inverse and we put $$r=\#\circ i^*$$ where
$i^*:{\cal G}^*\too S^*$ is the dual of the inclusion $i:S\inj\cal
G$.

With this observation in mind, the following proposition gives
another description of the solutions of the Yang-Baxter equation.
\begin{pr} Let $r\in {\cal G}\wedge{\cal G}$ and $(S_r,\om_r)$ its
associated subspace. The following assertions are equivalent:

1) $[r,r]=0.$

2) $r([\al,\be]_r)=[r(\al),r(\be)].$ $([\;,\;]_r$ is the bracket
given by  $(13))$.

3) $S_r$ is a subalgebra of $\cal G$ and $$\de\om_r(u,v,w):=\om_r(
u,[v,w])+\om_r(v ,[w,u])+\om_r( w,[u,v])=0$$for any $u,v,w\in
S_r$.

Moreover, if 1), 2) or 3) holds then $(\G^*,[\;,\;]_r)$ is a Lie
algebra and $r:\G^*\too\G$ is a morphism of Lie algebras.\end{pr}
 {\bf Proof:} The proposition follows from the following  formulas:
 $$\ga(r([\al,\be]_r)-[r(\al),r(\be)])=-[r,r](\al,\be,\ga),
 \qquad\al,\be,\ga\in\G^*\eqno(14)$$and, if
 $S$ is a subalgebra,
$$[r,r](\al,\be,\ga)=-\de\om_r(r(\al),r(\be),r(\ga)).\eqno(15)$$$
\Box$

This proposition shows that to give a solution of the Yang-Baxter
equation is equivalent to give a symplectic subalgebra of $\G$. We
recall that a symplectic algebra (see [3]) is a Lie   algebra $S$
endowed with a non-degenerate 2-form $\om$ such that $\de\om=0.$
\begin{Rm} Let $G$ be a Lie group,  $\G$ its  Lie algebra and
$S$ an even dimensional abelian subalgebra of $\G$. Any
non-degenerate 2-form $\om$ on $S$ verifies the assertion 3) in
Proposition 2.1 and hence $(S,\om)$ defines a solution of the
Yang-Baxter equation and then a structure of Poisson Lie group on
$G$.$\Box$\end{Rm}

The following proposition will be crucial in the simplification of
the equation $(11)$.

\begin{pr}
Let  $(\G,[\;\;],<,>)$ be a  Lie algebra endowed with a scalar
product,  $r\in\G\wedge\G$   a solution of $(Y-B)$ and
$(S_r,\om_r)$ its associated symplectic Lie algebra. Then,
$S_r\subset S_{<,>}$ iff the infinitesimal  Levi-Civita connection
$A^*$ associated with $(\G^*,[\;,\;]_r,<,>^*)$ is given by
$$A_\al^*\be=-ad^*_{r(\al)}\be,\qquad\al,\be\in\G^*.\eqno(16)$$

Moreover,   if $S_r\subset S_{<,>}$,  the curvature of $A^*$
vanishes and hence $(\G^*,[\;,\;]_r,<,>^*)$ is a Riemann-Lie
algebra.\end{pr}

{\bf Proof:}  $A^*$ is the unique bilinear map from $\G^*\times
\G^*$ to $\G^*$ such that:

1) $A^*_\al\be-A^*_\be\al=[\al,\be]_r$ for any $\al,\be\in\G^*$;

2) the endomorphism $A^*_\al:\G^*\too\G^*$ is skew-adjoint with
respect to  $<,>^*.$

The bilinear map $(\al,\be)\mapsto -ad^*_{r(\al)}\be$ verifies 1)
obviously and verifies 2)  iff $S_r\subset S_{<,>}$.

 If $A_\al^*\be=-ad^*_{r(\al)}\be$, the curvature of $A^*$ is
given by
$$
R(\al,\be)\ga=A^*_{[\al,\be]_r}\ga-
\left(A_\al^*A^*_\be\ga-A_\be^*A_\al^*\ga\right)
=ad_{r([\al,\be]_r)-[r(\al),r(\be)]}^*\ga=0$$from  Proposition 2.1
2). We conclude by Theorem 2.1.$\Box$
\begin{pr} Let $(\G,(\;,\;],<,>)$ be a Lie algebra with a scalar
 product. Let $r\in\G\wedge\G$ be a solution of $(Y-B)$ such that
 $S_r$ is a subalgebra of the orthogonal subalgebra $S_{<,>}$.
 Then $S_r$ is abelian.\end{pr}
 {\bf Proof:}
 $S_r$ is unimodular and symplectic and then resolvable (see
 [3]). Also $S_r$ carries a  bi-invariant scalar product so $S_r$
 must be abelian (see [5]).$\Box$

 We can now simplify the equation (11) and give the construction of
Riemann-Poisson Lie groups announced before.

Let $G$ be a Lie group,   $(\G,[\;,\;])$ its Lie algebra and $<,>$
 a scalar product  on $\G$.
We assume that the orthogonal subalgebra $S_{<,>}$ contains an
abelian even dimensional subalgebra $S$ endowed with   a
non-degenerate 2-form $\om$.

As in Remark 2.2, $(S,\om)$ defines a solution $r$ of $(Y-B)$ and
then a Poisson Lie tensor  $\pi$ on $G$. We can verify obviously
that, for any $g\in G$,
$$\pi^l(g)=r-Ad_g(r).$$
This implies that (11) can be rewritten
\begin{eqnarray*}
[Ad^*_g\left(A^*_\al\ga+ad^*_{r(\al)}\ga\right),Ad^*_g(\be)]_r
+[Ad^*_g(\al), Ad^*_g\left(A^*_\be\ga+ad^*_{r(\be)}\ga\right)]_r=\\
\;[Ad^*_g\left(ad^*_{Ad_g(r)(\al)}\ga\right),Ad^*_g(\be)]_r
+[Ad^*_g(\al), Ad^*_g\left(ad^*_{Ad_g(r)(\be)}\ga\right)]_r.
\end{eqnarray*}

Now, since $S\subset S_{<,>}$, we have by  Proposition 2.2
$$A^*_\al\ga+ad^*_{r(\al)}\ga=0$$for any
$\al,\ga\in\G^*$. On other hand, the following formula is easy to
get
$$Ad^*_g[ad_{r(\al)}^*\be]=
ad_{(Ad_{g^{-1}}r)(Ad^*_g\al)}^*(Ad_g^*\be),\qquad g\in
G,\al,\be\in\G^*.$$ Finally, $(G,\pi,<,>^{*l})$ is a
Riemann-Poisson Lie group iff
$$[ad_{r(\al)}^*\ga,\be]_r+[\al,ad_{r(\be)}^*\ga]_r=0,\qquad\al,\be,\ga\in\G^*.$$
But, also since $A^*_\al\ga+ad^*_{r(\al)}\ga=0,$ this condition is
equivalent to  $(\G^*,[\;\;]_r,<,>^*)$ is a Riemann-Lie algebra
which is true by Proposition 2.2. So, we have shown:
\begin{th} Let $G$ be a Lie group,  $(\G,[\;,\;])$ its Lie
algebra and $<,>$ a scalar product on $\G$. Let $S$ be an even
dimensional abelian subalgebra of the orthogonal subalgebra
$S_{<,>}$ and $\om$  a non-degenerate 2-form on $S$. Then, the
solution of the Yang-Baxter equation associated with $(S,\om)$
defines a structure of Poisson Lie group $(G,\pi)$ and
$(G,\pi,<,>^{*l})$ is a Riemann-Poisson Lie group.\end{th}

Let us enumerate some important  cases where  this theorem can be
used.

 1) Let $G$ be a compact Lie group and $\G$ its Lie algebra.
 For any
  bi-invariant scalar
product $<,>$ on the Lie algebra $\G$, $S_{<,>}=\G$. By Theorem
2.2, we can associate to any even dimensional subalgebra of $\G$ a
Riemann-Poisson Lie group structure on $G$.

2) Let $(\G,[\;,\;],<,>)$ be a Riemann-Lie algebra. By Theorem
2.1, the orthogonal subalgebra $S_{<,>}$ is abelian and any even
dimensional subalgebra of $S_{<,>}$ gives rise to a structure of a
Riemann-Poisson Lie group on any Lie group whose the Lie algebra
is $\G$. Moreover, we get a structure of  bialgebra
$(\G,[\;,\;],\G^*,[\;,\;]_r)$ where both $\G$ and $\G^*$ are
Riemann-Lie algebras.

Finally, we observe that the Riemann-Lie groups constructed above
inherit the properties of Riemann-Poisson manifolds (see [2]).
Namely, the symplectic leaves of these Poisson Lie groups are
K\"ahlerian and their Poisson structures  are unimodular.

\section{Proof of the equivalence ``$1)\Leftrightarrow5)$''
in Theorem 2.1}

In this section we will give a proof of the equivalence
``$1)\Leftrightarrow5)$'' in Theorem 2.1. The proof is a sequence
of lemmas. Namely,  we will show that, for a Riemann-Lie algebra
$(\G,[\;,\;],<,>)$, the orthogonal subalgebra $S_{<,>}$ is
abelian. Moreover, $S_{<,>}$ is the $<,>$-orthogonal of the ideal
$[\G,\G]$. This result will be the key of the proof.

We begin by  a characterization of Riemann-Lie subalgebras.

\begin{pr} Let $(\G,[\;,\;],<,>)$ be a Riemann-Lie algebra and
$\cal H$ a subalgebra of $\G$. For any $u,v\in\cal H$, we put
$A_uv=A_u^0v+A_u^1v$, where $A_u^0v\in\cal H$ and $A_u^1v\in{\cal
H}^\perp$. Then, $({\cal H},[\;,\;],<,>)$ is a Riemann-Lie algebra
if and only if, for any $u,v,w\in\cal H$,
$[A^1_uv,w]+[v,A^1_uw]\in{\cal H}^\perp.$\end{pr}

{\bf Proof:}  We have, from (9), that for any $u,v,w\in\cal H$
$$[u,[v,w]]=[A^0_uv,w]+[v,A^0_uw]+[A^1_uv,w]+[v,A^1_uw].$$Now
$A^0:{\cal H}\times{\cal H}\too\cal H$ is the infinitesimal
Levi-Civita connection associated with the restriction of $<,>$ to
$\cal H$ and the proposition follows.$\Box$

 We will introduce now some mathematical objects which will be useful
latter.

 Let $(\G,[\;,\;],<,>)$ a Lie algebra
endowed with a scalar product.

From (5), we  deduce obviously that the infinitesimal Levi-Civita
connection $A$ associated to $<,>$ is given by

$$A_uv=\frac12[u,v]-\frac12\left(ad_u^tv+ad_v^tu\right)\quad
u,v\in\G.\eqno(17)$$

On other hand, the orthogonal with respect to $<,>$ of the ideal
$[\G,\G]$ is given by
$$[\G,\G]^\perp=\bigcap_{u\in\G}Kerad^t_u.\eqno(18)$$

Let us introduce, for any $u\in\G$, the endomorphism
$$D_u=ad_u-A_u.\eqno(19)$$ We have, by a straightforward calculation,
the relations
\begin{eqnarray*}
D_u(v)&=&\frac12[u,v]+\frac12\left(ad_u^tv+ad_v^tu\right),\\
D_u^t(v)&=&\frac12[u,v]+\frac12\left(ad_u^tv-ad_v^tu\right).
\end{eqnarray*}
From these relations, we remark that, for any $u,v\in\G$,
$D_u^t(v)=-D_v^t(u)$ and then $$\forall
u\in\G,\;D_u^t(u)=0.\eqno(20)$$ We remark also that
$$D_u^t=D_u\quad\Leftrightarrow\quad\forall v\in\G,
ad_v^tu=0.$$So, by (18), we get
$$[\G,\G]^\perp=\{u\in\G; D_u^t=D_u\}.\eqno(21)$$

We begin now to state and to prove a sequence of results which
will give a proof of the equivalence ``$1)\Leftrightarrow5)$'' in
Theorem 2.1.
\begin{pr} Let $(\G,[\;,\;],<,>)$ be a Riemann-Lie algebra. Then
$Z(\G)^\perp$ ($Z(\G)$ is the center of $\G$) is an ideal of $\G$
which contains the ideal $[\G,\G]$. In particular,
$$\G=Z(\G)\oplus Z(\G)^\perp.$$\end{pr}

{\bf Proof:} For any $u\in Z(\G)$ and $v\in\G$,  from (17) and the
fact that $A_u$ is skew-adjoint, $A_uv=-\frac12ad_v^tu\in
Z(\G)^\perp.$ By (8),  for any $w\in\G$
$$[A_uv,w]=[A_wv,u]=0,$$so $A_uv\in Z(\G)$ and then
$A_uv=-\frac12ad_v^tu=0$ which shows that $u\in[\G,\G]^\perp$. So
 $Z(\G)\subset[\G,\G]^\perp$ and the proposition
follows.$\Box$

From this proposition and the fact that for a nilpotent Lie
algebra $\G$  $Z(\G)\cap[\G,\G]\not=\{0\}$, we get the following
lemma.
\begin{Le} A nilpotent Lie algebra $\G$ carries a structure of
Riemann-Lie algebra if and only if $\G$ is abelian.\end{Le}

We can now get the following crucial result.
\begin{Le} Let $(\G,[\;,\;],<,>)$ be a Riemann-Lie algebra. Then
the orthogonal Lie subalgebra $S_{<,>}$ is  abelian.
\end{Le}

{\bf Proof:} By (17),  $A_uv=\frac12[u,v]$ for any $u,v\in
S_{<,>}$. So,  by Proposition 3.1, $S_{<,>}$ is a
Riemann-subalgebra. By (9), we have, for any $u,v,w\in S_{<,>}$,
\begin{eqnarray*}
[u,[v,w]]&=&[A_uv,w]+[v,A_uw]\\
&=&\frac12[[u,v],w]+\frac12[v,[u,w]]\\
&=&\frac12[u,[v,w]]\end{eqnarray*}and then
$[S_{<,>},[S_{<,>},S_{<,>}]]=0$ i.e. $S_{<,>}$ is a nilpotent Lie
algebra and then abelian by Lemma 3.1.$\Box$

\begin{Le} Let $(\G,[\;,\;],<,>)$ be a Riemann-Lie algebra. Then
$$[\G,\G]^\perp=\{u\in\G; D_u=0\}.$$\end{Le}

{\bf Proof:} Firstly, we notice that, by $(21)$,
$[\G,\G]^\perp\supset\{u\in\G; D_u=0\}.$ On other hand,  remark
that the relation (8) can be rewritten
$$[D_u(v),w]+[v,D_u(w)]=0$$for any $u,v,w\in\G$. So, we can  deduce
immediately that $[Ker D_u,ImD_u]=0$ for any $u\in\G$.

Now we observe that, for any $u\in[\G,\G]^\perp$, the endomorphism
$D_u$ is auto-adjoint and then diagonalizeable on $\reel$. Let
$u\in[\G,\G]^\perp$,  $\la\in\reel$ be an eigenvalue of $D_u$ and
$v\in\G$ an eigenvector associated with $\la$. We have
$$<D_u(v),v>=\la<v,v>\stackrel{(\al)}{=}-<A_vu,v>\stackrel{(\be)}{=}
-<[v,u],v>\stackrel{(\ga)}{=}0.$$So $\la=0$ and we obtain that
$D_u$ vanishes identically. Hence the lemma follows.

The equality $(\al)$ is a consequence of the definition of $D_u$,
and the equality  $(\be)$ follows from the definition of $A$. We
observe that $v\in ImD_u$ and $u\in KerD_u$  since
$D_u(u)=D^t_u(u)=0$ (see (20)) and the equality $(\ga)$ follows
from the remark above. $\Box$
\begin{Le} Let $(\G,[\;,\;],<,>)$ be a Riemann-Lie algebra. Then
$$S_{<,>}=[\G,\G]^\perp.$$\end{Le}

{\bf Proof:} From Lemma 3.3,  for any $u\in[\G,\G]^\perp$,
$A_u=ad_u$ and then $ad_u$ is skew-adjoint.  So
$[\G,\G]^\perp\subset S_{<,>}$. To prove the second inclusion we
need to work harder than the first one.

Firstly, remark that one can suppose that $Z(\G)=\{0\}$. Indeed,
$\G=Z(\G)\oplus Z(\G)^\perp$ (see Proposition 3.2), $Z(\G)^\perp$
is a Riemann-Lie algebra (see Proposition 3.1),
$[\G,\G]=[Z(\G)^\perp,Z(\G)^\perp]$ and $S_{<,>}=Z(\G)\oplus
S'_{<,>}$ where $S'_{<,>}$ is the orthogonal subalgebra associated
to $(Z(\G)^\perp,<,>)$.

We suppose now that $(\G,[\;,\;],<,>)$ is a Riemann-Lie algebra
such that $Z(\G)=\{0\}$ and we want   to prove the inclusion
$[\G,\G]^\perp\supset S_{<,>}$. Notice that it suffices to show
 that, for any $u\in S_{<,>}$, $A_u=ad_u$.

 The proof requires some preparation. Let us introduce
 the subalgebra $K$ given
by
$$K=\bigcap_{u\in S_{<,>}}\ker ad_u.$$

Firstly, we notice that $K$ contains $S_{<,>}$ because $S_{<,>}$
is abelian (see Lemma 3.2).

 On other hand, we remark that, for any
$u\in S_{<,>}$, the endomorphism $A_u$ leaves invariant $K$ and
$K^\perp$.  Indeed, for any $v\in K$ and any $w\in S_{<,>}$, we
have
$$
[w,A_uv]\stackrel{(\al)}{=}[w,A_vu]\stackrel{(\be)}{=}-[A_wu,v]
\stackrel{(\ga)}{=}0$$ and then $A_uv\in K$, this shows that $A_u$
leaves invariant $K$. Furthermore, $A_u$ being skew-adjoint, we
have $A_u(K^\perp)\subset K^\perp$.

The equality $(\al)$ follows from the relation
$A_uv=A_vu+[u,v]=A_uv$, the equality $(\be)$ follows from (8) and
$(\ga)$ follows from the relation $A_wu=\frac12[w,u]=0.$

With this observation in mind, we consider now the representation
$\rho:S_{<,>}\too so(K^\perp)$ given by
$$\rho(u)={ad_u}_{| K^\perp}\qquad u\in S_{<,>}.$$
It is clear that
$$\bigcap_{u\in S_{<,>}}\ker\rho(u)=\{0\}.\eqno(*)$$This relation and the
fact that $S_{<,>}$ is abelian imply that $\dim K^\perp$ is even
and that there is an orthonormal basis $(e_1,f_1,\ldots,e_p,f_p)$
of $K^\perp$ such that
$$\forall i\in\{1,\ldots, p\},\forall u\in S_{<,>},\quad
ad_u{e_i}=\la^i(u)f_i\quad\mbox{and}\quad
ad_uf_i=-\la^i(u)e_i,\eqno(**)$$where $\la^i\in S_{<,>}^*$.

Now, for any $u\in S_{<,>}$, since $A_u$ leaves $K^\perp$
invariant, we can write
\begin{eqnarray*}
A_ue_i&=&\sum_{j=1}^p\left(<A_ue_i,e_j>e_j+<A_ue_i,f_j>f_j\right),\\
A_uf_i&=&\sum_{j=1}^p\left(<A_uf_i,e_j>e_j+<A_uf_i,f_j>f_j\right).
\end{eqnarray*}
From (9), we have  for any $v\in S_{<,>}$ and for any
$i\in\{1,\ldots, p\}$,
\begin{eqnarray*}
\;[u,[v,e_i]]&=&[A_uv,e_i]+[v,A_ue_i],\\
\;[u,[v,f_i]]&=&[A_uv,f_i]+[v,A_uf_i].\end{eqnarray*} Using the
the equality  $A_uv=0$ and $(**)$ and substituting we get
\begin{eqnarray*}
-\la^i(u)\la^i(v)e_i&=&\sum_{j=1}^p\la^j(v)<A_ue_i,e_j>f_j-
\sum_{j=1}^p\la^j(v)<A_ue_i,f_j>e_j,\\
-\la^i(u)\la^i(v)f_i&=&\sum_{j=1}^p\la^j(v)<A_uf_i,e_j>f_j-
\sum_{j=1}^p\la^j(v)<A_uf_i,f_j>e_j.\end{eqnarray*} Now, it is
clear from $(*)$  that, for any $i\in\{1,\ldots, p\}$, there
exists $v\in S_{<,>}$ such that $\la^i(v)\not=0$. Using this fact
and the relations above, we get
$$A_ue_i=\la^i(u)f_i\quad\mbox{and}\quad A_uf_i=-\la^i(u)e_i.$$So
we have shown that, for any $u\in S_{<,>}$,
$${A_u}_{| K^\perp}={ad_u}_{| K^\perp}.$$

 Now, for any $u\in S_{<,>}$ and for any $k\in K$, $ad_u(k)=0$.
 So,
 to complete the proof of the lemma, we will show that, for any
$u\in S_{<,>}$ and for any $k\in K$,  $A_uk=0.$ This will be done
by showing that $A_uk\in Z(\G)$ and conclude by using the
assumption $Z(\G)=\{0\}.$

Indeed, for any $h\in K$, by (8)
$$[A_uk,h]=[A_hk,u].$$  Since $A_u(K)\subset
K$ and since $K$ is a subalgebra, $[A_uk,h]\in K$. Now, $K\subset
\ker ad_u$ and $ad_u$ is skew-adjoint so $[A_hk,u]\in Im
ad_u\subset K^\perp$  . So $[A_uk,h]=0.$ On other hand, for any
$f\in K^\perp$, we have, also from (8),
$$[A_uk,f]=[A_ku,f]=[A_fu,k]=0$$since
$A_fu=[f,u]+A_uf=[f,u]+[u,f]=0.$

We deduce that $A_uk\in Z(\G)$ and then $A_uk=0$. The proof of the
lemma is complete.$\Box$
\begin{Le} Let $(\G,[\;,\;],<,>)$ be a Riemann-Lie algebra such
that $Z(\G)=0.$ Then
$$\G\not=[\G,\G].$$\end{Le}

{\bf Proof:} Let $(\G,[\;,\;],<,>)$ be a Riemann-Lie algebra such
that $Z(\G)=0.$ We will show that the assumption $\G=[\G,\G]$
implies that the Killing form of $\G$ is strictly negative
definite  and then $\G$ is semi-simple and compact which is in
contradiction with  lemma 2.1.

Let $u\in\G$ fixed. Since $A_u$ is skew-adjoint, there is an
orthonormal basis $(a_1,b_1,\ldots,a_r,b_r,c_1,\ldots, c_l)$ of
$\G$ and $(\mu_1,\ldots,\mu_r)\in\reel^r$ such that, for any $
i\in\{1,\ldots,r\}$ and any $j\in\{1,\ldots,l\}$,
$$
A_ua_i=\mu_ib_i,\quad A_ub_i=-\mu_ia_i\quad\mbox{and}\quad
A_uc_j=0.$$Moreover, $\mu_i>0$ for any $i\in\{1,\ldots,r\}$.

By applying $(9)$, we can deduce, for any $i,j\in\{1,\ldots,r\}$
and for any $k,h\in\{1,\ldots, l\}$, the relations:
\begin{eqnarray*}
\;[u,[a_i,a_j]]&=&\mu_i[b_i,a_j]+\mu_j[a_i,b_j],\quad
\;[u,[b_i,b_j]]=-\mu_j[b_i,a_j]-\mu_i[a_i,b_j],\\
\;[u,[a_i,b_j]]&=&-\mu_j[a_i,a_j]+\mu_i[b_i,b_j],\quad
\;[u,[b_i,a_j]]=-\mu_i[a_i,a_j]+\mu_j[b_i,b_j],\\
\;[u,[c_k,a_j]]&=&\mu_j[c_k,b_j],\quad
\;[u,[c_k,b_j]]=-\mu_j[c_k,a_j],\quad
\;[u,[c_k,c_h]]=0.\end{eqnarray*} From these relations we deduce
\begin{eqnarray*}
ad_u\circ ad_u([a_i,a_j])&=&
-(\mu_i^2+\mu_j^2)[a_i,a_j]+2\mu_i\mu_j[b_i,b_j],\\ ad_u\circ
ad_u([b_i,b_j])&=&
2\mu_i\mu_j[a_i,a_j]-(\mu_i^2+\mu_j^2)[b_i,b_j],\\
ad_u\circ ad_u([b_i,a_j])&=&
-(\mu_i^2+\mu_j^2)[b_i,a_j]-2\mu_i\mu_j[a_i,b_j],\\ ad_u\circ
ad_u([a_i,b_j])&=&
-2\mu_i\mu_j[b_i,a_j]-(\mu_i^2+\mu_j^2)[a_i,b_j],\\
ad_u\circ ad_u([c_k,a_j])&=&-\mu_j^2[c_k,a_j],\\ ad_u\circ
ad_u([c_k,b_j])&=&-\mu_j^2[c_k,b_j],\\ ad_u\circ
ad_u([c_k,c_h])&=&0.\end{eqnarray*} By an obvious transformation
we obtain
\begin{eqnarray*}
ad_u\circ ad_u\left([a_i,a_j]+[b_i,b_j]\right)&=&-(\mu_i-\mu_j)^2
\left([a_i,a_j]+[b_i,b_j]\right),\\
ad_u\circ ad_u\left([a_i,a_j]-[b_i,b_j]\right)&=&-(\mu_i+\mu_j)^2
\left([a_i,a_j]-[b_i,b_j]\right),\\
ad_u\circ ad_u\left([b_i,a_j]+[a_i,b_j]\right)&=&-(\mu_i+\mu_j)^2
\left([b_i,a_j]+[a_i,b_j]\right),\\
ad_u\circ ad_u\left([b_i,a_j]-[a_i,b_j]\right)&=&-(\mu_i-\mu_j)^2
\left([b_i,a_j]-[a_i,b_j]\right),\\
ad_u\circ ad_u([c_k,a_j])&=&-\mu_j^2[c_k,a_j],\\ ad_u\circ
ad_u([c_k,b_j])&=&-\mu_j^2[c_k,b_j],\\ ad_u\circ
ad_u([c_k,c_h])&=&0.\end{eqnarray*}

Suppose now  $\G=[\G,\G]$. Then the family of vectors

$\left\{[a_i,a_j]+[b_i,b_j],[a_i,a_j]-[b_i,b_j],[b_i,a_j]+
[a_i,b_j],\right.$

$[b_i,a_j]- [a_i,b_j], [c_k,a_i],[c_k,b_j],[c_k,c_h];$
$\left.i,j\in\{1,\ldots,r\},h,k\in\{1,\ldots,l\}\right\}$ spans
$\G$ and then $ad_u\circ ad_u$ is diagonalizeable and all its
eigenvalues are non positive. Now its easy to deduce that
$ad_u\circ ad_u=0$ if and only if $ad_u=0$. Since $Z(\G)=0$ we
have shown that, for any $u\in \G\setminus\{0\}$, $Tr(ad_u\circ
ad_u)<0$ and then the Killing form of $\G$ is strictly negative
definite and then $\G$ is semi-simple compact. We can conclude
with the Lemma 2.1.$\Box$

{\bf Proof of the equivalence ``$1)\Leftrightarrow5)$'' in Theorem
2.1.}

It is an obvious and straightforward calculation to show that
$5)\Rightarrow1)$.

Conversely, let $(\G,[\;,\;],<,>)$ be a Riemann-Lie algebra. By
Proposition 3.2, we can suppose that $Z(\G)=\{0\}$.

We have, from Lemma 3.5 and Lemma 3.4, $\G\not=[\G,\G]$ which
implies $S_{<,>}\not=0$ and
$\G=S_{<,>}\stackrel{\perp}\oplus[\G,\G]$. Moreover, $[\G,\G]$ is
a Riemann-Lie algebra (see Proposition 3.1) and we can repeat the
argument above to deduce that finally $\G$ is solvable which
implies that $[\G,\G]$ is nilpotent and then abelian by Lemma 3.1
and the implication follows.$\Box$

\begin{Rm}  The pseudo-Riemann-Lie algebras are completely
different from the Riemann-Lie algebras. Indeed,  the
3-dimensional Heisenberg Lie algebra which is nilpotent carries a
Lorentzian-Lie algebra structure. On other hand, the non trivial
2-dimensional Lie algebra carries a Lorentzian inner product whose
curvature vanishes and don't carries any structure of
pseudo-Riemann-Lie algebra.$\Box$
\end{Rm}\bigskip

{\bf References}

\bigskip

[1] M. Boucetta,  Compatibilit\'e des structures
pseudo-riemanniennes et des structures de Poisson, C. R. Acad.
Sci. Paris, {\bf t. 333}, S\'erie I, (2001) 763-768.

[2] M. Boucetta, Poisson manifolds with compatible pseudo-metric
and pseudo-Riemannian Lie algebras, Preprint math.DG/0206102. To
appear in  Differential Geometry and its Applications.

[3] A. Lichnerowicz and A. Medina, On Lie groups with
left-invariant symplectic or K\"ahlerian structures, Letters in
Mathematical Physics {\bf16} (1988), 225-235.

[4] A. Medina, P. Revoy,  Alg\`ebres de Lie et produit scalaire
invariant, Ann. Ec. Norm. Sup., 4\`eme s\'erie, {\bf t. 18} (1985)
553-561.

[5] J. Milnor, Curvature of Left invariant metrics on Lie groups,
Advances in Mahtematics {\bf21}, 293-329 (1976).

[6] J. H. Lu and A. Weinstein, Poisson Lie groups, dressing
transformations and Bruhat decompositions, J. Differential
Geometry, {\bf31} (1990) 501-526.

\bigskip

{\it Mohamed BOUCETTA

Facult\'e des sciences et techniques, BP 549, Marrakech, Maroc

Email: boucetta@fstg-marrakech.ac.ma}

\end{document}